\documentclass[10pt]{amsart}
\usepackage{amsmath,amstext,amssymb,amsopn,amsthm,mathrsfs}

\newtheorem{theorem}{Theorem}[section]

\newtheorem{corollary}{Corollary}[section]
\newtheorem{lemma}{Lemma}[section]
\newtheorem{proposition}{Proposition}[section]
\newtheorem{obs}{Observation}[section]
\numberwithin{equation}{section}
\newcommand{\dem}{\medskip \par \noindent \mbox{\bf Proof. }}
\def\ep{\hfill{$\Box $}}

\begin{document}

\title[Jacobi polynomials and the Watson Kernel]{On Abel summability of Jacobi polynomials series, the Watson Kernel and applications}

\author{Calixto P. Calder\'on}
\address{Department of Mathematics, University of Illinois at Chicago, Il, 60607, USA.}
\email{[Calixto P. Calder\'on]cpc@uic.edu}
\author{Wilfredo O. Urbina}
\address{
 Department of Mathematical and Actuarial Sciences, Roosevelt University  Chicago, Il, 60605, USA.}
\email{[Wilfredo Urbina]wurbinaromero@roosevelt.edu}

\begin{abstract}
In this paper we return to the study of the Watson kernel for the Abel summabilty of Jacobi polynomial series. These estimates have been studied for over more than 30 years. The main innovations are in the techniques used to get the estimates that allow us to handle the case $0<\alpha $ as well as $ -1< \alpha <0$, with essentially the same method; using an integral superposition of Poisson type kernel and Muckenhoupt $A_p$-weight theory. We consider a generalization of a theorem due to Zygmund in the context to Borel measures. The proofs are therefore different from the ones given in \cite{cafacal}, \cite{cafacal1}, \cite{cafacal2}  and  \cite{calv}. We will also discuss in detail the Calder\'on-Zygmund decomposition for non-atomic Borel measures in $\mathbb{R}$. Then, we prove that the Jacobi measure is doubling and therefore, following \cite{apCal}, we study the corresponding $A_p$  weight theory in the setting of Jacobi expansions, considering  power weights of the form $(1-x)^{\overline{\alpha}}, \, (1+x)^{\overline{\beta}}$, $-1 < {\overline{\alpha}}<0,\,  -1 < {\overline{\beta}}<0$ with negative exponents. Finally, as an application of the weight theory we obtain $L^p$ estimates for the maximal operator of Abel summability of Jacobi function expansions for suitable values of $p$.
\end{abstract}
\maketitle

\section{Introduction}
Given $\alpha, \beta >-1,$  consider the Jacobi measure $J^{\alpha,\beta}$ on $[-1,1]$, defined as
\begin{equation}
J^{\alpha,\beta} (dx) =\omega_{\alpha,\beta}(x)dx = (1-x)^{\alpha}(1+x)^{\beta} dx.
\end{equation}
 The Jacobi polynomials of parameters $\alpha, \beta$, $\{P^{\alpha,\beta}_n\}_{n\geq0}$ are the orthogonal polynomials with respect to the measure $J^{\alpha,\beta},$
\begin{eqnarray}\label{orthopoly}
\int_{-1}^1 P^{\alpha,\beta}_n(x)  P^{\alpha,\beta}_m(x) (1-x)^{\alpha}(1+x)^{\beta} dx = 0,\quad \mbox{if} \quad n\neq m,
\end{eqnarray}
with
\begin{eqnarray}
 \nonumber \int_{-1}^1 [P^{\alpha,\beta}_n(x)]^2 (1-x)^{\alpha}(1+x)^{\beta} dx &=&\frac{2^{\alpha+\beta+1}}{2n+\alpha+\beta+1} \frac{\Gamma(n+\alpha+1)\Gamma(n+\beta+1)}{\Gamma(n+1)\Gamma(n+\alpha\beta+1)} \\
&=& h^{(\alpha,\beta)}_n.
\end{eqnarray}

The normalization is given by 
\begin{equation}
 P^{\alpha,\beta}_n(1) = {n+\alpha \choose n}
\end{equation}

The Jacobi functions are defined, for each $n$, as 
\begin{equation}\label{JacobiFunct}
F^{(\alpha,\beta)}_n(x) = P^{\alpha,\beta}_n(x) (1-x)^{\alpha/2}(1+x)^{\beta/2};
\end{equation}
therefore, from (\ref{orthopoly}) one gets that the Jacobi functions $\{F^{(\alpha,\beta)}_n\}$  are orthogonal on       $[-1,1]$  with respect to the Lebesgue measure,
\begin{eqnarray*}
\int_{-1}^1 F^{\alpha,\beta}_n(x)  F^{\alpha,\beta}_m(x)  dx = 0,\quad \mbox{if} \quad n\neq m.
\end{eqnarray*}
For any $f \in L^2([-1,1], J^{\alpha,\beta})$ we consider its Fourier-Jacobi polynomial expansion
\begin{equation}\label{FJexpP}
\sum_{n=0}^{\infty} \hat{f}^{(\alpha,\beta)}(n) P^{\alpha,\beta}_n(x),
\end{equation}
where
$$  \hat{f}^{(\alpha,\beta)}(n) = \frac{1}{h^{(\alpha,\beta)}_n} \int_{-1}^1f(y)  P^{\alpha,\beta}_n(y)  J^{\alpha,\beta}(dy),$$
is the $n$-th Fourier-Jacobi polynomial  coefficient. Then its partial sum $s_m^{\alpha,\beta}(f,x)$, can be written as
\begin{equation}
s_m^{\alpha,\beta}(f,x) =  \int_{-1}^{1} {\mathcal K}_m^{\alpha,\beta}(x,y) f(y) J^{\alpha,\beta}(dy)
\end{equation}
where 
\begin{eqnarray*}
  {\mathcal K}_m^{\alpha,\beta}(x,y)&=&\sum_{n=0}^{m}  \frac{ P^{\alpha,\beta}_n(x) P^{\alpha,\beta}_n(y)}{h^{(\alpha,\beta)}_n}.
\end{eqnarray*}

The kernel $ {\mathcal K}_m^{\alpha,\beta}$ is called  the {\it Dirichlet-Jacobi kernel.}

By using the three term recurrent relation of $\{P^{\alpha,\beta}_n\}$ it is well know that 
one can get the Christofel-Darboux formula for $ {\mathcal K}_m^{\alpha, \beta},$
$$ {\mathcal K}_m^{\alpha, \beta}(x,y) = \frac{2^{-\alpha-\beta}}{2m+\alpha+\beta+2} \frac{\Gamma(m+2) \Gamma(m+\alpha+\beta+2)}{\Gamma(m+\alpha+1)\Gamma(m+\beta+1)} \frac{P_{m+1}^{\alpha, \beta}(x) P_m^{\alpha, \beta}(y) -P_m^{\alpha, \beta}(x) P_{m+1}^{\alpha, \beta}(y) }{x-y}.$$
Moreover by orthogonality, we get
$$\int_{-1}^1  {\mathcal K}_m^{\alpha, \beta}(x,y) J^{\alpha,\beta}(dy) =1.$$

Let us consider the Abel summability of the Jacobi polynomial series expansion of $f$ (\ref{FJexpP}),
\begin{equation}\label{abel}
f^{\alpha, \beta}(r,x) = \sum_{n=0}^{\infty}r^n \hat{f}^{(\alpha,\beta)}(n) P^{\alpha,\beta}_n(x), \quad 0 < r < 1.
\end{equation}
Using  a classical argument and the estimate (see \cite{sz} (7.32.1) )
\begin{equation}\label{Jacobibound}
|P^{\alpha,\beta}_n(x)| \leq C n^{q+1/2},
\end{equation}
where $q= \max(\alpha, \beta) \geq -1/2$,   it is easy to see that the series (\ref{abel}) converges uniformly and absolutely in $[-1,1]$; and  therefore $f^{\alpha, \beta}(r,x)$ has an integral representation,
\begin{equation}
f^{\alpha, \beta}(r,x) = \int_{-1}^{1} K^{\alpha,\beta}(r,x,y) f(y) J^{\alpha,\beta}(dy),
\end{equation}
for $f \in L^1([-1,1], J_{\alpha,\beta})$, where 
\begin{eqnarray}
K^{\alpha,\beta}(r,x,y) &=& \sum_{n=0}^{\infty} r^n \frac{ P^{\alpha,\beta}_n(x) P^{\alpha,\beta}_n(y)}{h^{(\alpha,\beta)}_n}.
\end{eqnarray}
$K^{\alpha,\beta}$ is called the {\em Watson kernel}. Observe that trivially the kernel is symmetric in $x$ and $y$, i.e.
$ K^{(\alpha,\beta)}(r,x,y) =K^{(\alpha,\beta)}(r,y,x).$

Analogously, for any $f \in L^2([-1,1])$ we consider its Fourier-Jacobi function expansion
\begin{equation}\label{FJexpF}
\sum_{n=0}^{\infty} \tilde{f}^{(\alpha,\beta)}(n) F^{\alpha,\beta}_n(x),
\end{equation}
where
$$  \tilde{f}^{(\alpha,\beta)}(n) = \frac{1}{h^{(\alpha,\beta)}_n} \int_{-1}^1f(y)  F^{\alpha,\beta}_n(y)  dy,$$
is the $n$-th Fourier-Jacobi function coefficient. Then its partial sum $\tilde{s}_m^{\alpha,\beta}(f,x)$, can be written as
\begin{equation}
\tilde{s}_m^{\alpha,\beta}(f,x) =  \int_{-1}^{1} \tilde{{\mathcal K}}_m^{\alpha,\beta}(x,y) f(y) dy,
\end{equation}
where 
\begin{eqnarray*}
 \tilde{{\mathcal K}}_m^{\alpha,\beta}(x,y)&=&\sum_{n=0}^{m}  \frac{ F^{\alpha,\beta}_n(x) F^{\alpha,\beta}_n(y)}{h^{(\alpha,\beta)}_n}\\
 &=&  \sum_{n=0}^{m} \frac{ P^{\alpha,\beta}_n(x) P^{\alpha,\beta}_n(y)}{h^{(\alpha,\beta)}_n}(1-x)^{\alpha/2}(1-y)^{\alpha/2}(1+x)^{\beta/2}(1+y)^{\beta/2}.
\end{eqnarray*}
Now consider the Abel summability of the Jacobi function series expansion of $f$ (\ref{FJexpF}),
\begin{equation}
 \tilde{f}^{\alpha, \beta}(r,x) = \sum_{n=0}^{\infty}r^n \tilde{f}^{(\alpha,\beta)}(n) F^{\alpha,\beta}_n(x), \quad 0 < r < 1,
\end{equation}
then we also get an integral representation,
\begin{equation}\label{JabFuncIntRep}
 \tilde{f}^{\alpha, \beta}(r,x) = \int_{-1}^{1} \tilde{K}^{\alpha,\beta}(r,x,y) f(y) dy,
\end{equation}
for $f \in L^1([-1,1])$, where 
\begin{eqnarray*}
\tilde{K}^{\alpha,\beta}(r,x,y) &=& \sum_{n=0}^{\infty} r^n \frac{ F^{\alpha,\beta}_n(x) F^{\alpha,\beta}_n(y)}{h^{(\alpha,\beta)}_n}\\
&=&  \sum_{n=0}^{\infty} r^n \frac{ P^{\alpha,\beta}_n(x) P^{\alpha,\beta}_n(y)}{h^{(\alpha,\beta)}_n}(1-x)^{\alpha/2}(1-y)^{\alpha/2}(1+x)^{\beta/2}(1+y)^{\beta/2}\\
&=& K^{\alpha,\beta}(r,x,y) (1-x)^{\alpha/2}(1-y)^{\alpha/2}(1+x)^{\beta/2}(1+y)^{\beta/2}.
\end{eqnarray*}
$\tilde{K}^{\alpha,\beta}$ is called the {\em modified Watson kernel} for Jacobi functions.\\

From the previous representation and (\ref{JabFuncIntRep}) we get,
\begin{equation}\label{JabFuncIntRep2}
 \tilde{f}^{\alpha, \beta}(r,x) = (1-x)^{\alpha/2} (1+x)^{\beta/2} \int_{-1}^{1}K^{\alpha,\beta}(r,x,y) (1-y)^{\alpha/2}(1+y)^{\beta/2} f(y) dy.\\
\end{equation}

 In 1936 Watson obtained  the following representation for $K^{\alpha,\beta}(r,x,y)$, see \cite{bate} page 272,
\begin{equation}\label{Watsonrep}
K^{\alpha,\beta}(r,x,y) = r^{(1-\alpha-\beta)/2} \frac{d}{dr}\left( k^{1+\alpha+\beta} \int_0^{\pi/2} \frac{sec^{2+\alpha+\beta} \omega \, \cos (\alpha-\beta) \omega}{Z_1^\alpha Z_2^\beta Y} d\omega\right)
\end{equation}
where $ k = \frac{1}{2} (r^{1/2} + r^{-1/2}), s = k \sec \omega$,
$$ Y = \left(\left(\frac{x-y}{2}\right)^2+ (s^2-1)(s^2-xy) \right)^{1/2},$$
$$Z_1 = s^2-\frac{1}{2} (x+y) +Y, \, \mbox{and}$$
$$Z_2 = s^2+\frac{1}{2} (x+y) +Y.$$
The integral in (\ref{Watsonrep}) can be proved that is convergent only if $\alpha+\beta>-1;$ since $s \geq 2$, $Y^2 \sim s^4, Z_1 \sim s^2, Z_2 \sim s^2$, then taking the change of variable $s = k \sec \omega$,

$$\int_0^{\pi/2} \frac{sec^{2+\alpha+\beta} \omega \, \cos (\alpha-\beta) \omega}{Z_1^\alpha Z_2^\beta Y} d\omega 
\leq k^{-(2+\alpha+\beta)}  \int_k^{\infty} \frac{s^{\alpha+\beta+1}}{Z_1^{\alpha} Z_2^{\beta} Y}\frac{ k \, ds}{s \sqrt{s^2-k^2}},$$
Assuming that $ 1/2 < r < 1$, and then $1 < k <3/2 < 2,$ for $2 < s < \infty$, 
\begin{eqnarray}\label{estimate1}
 \int_2^{\infty} \frac{s^{\alpha+\beta+1}}{Z_1^{\alpha} Z_2^{\beta} Y}\frac{ k \, ds}{s \sqrt{s^2-k^2}}
&\sim& 
\nonumber C \int_k^{\infty} \frac{s^{\alpha+\beta+1}}{s^{2\alpha} s^{2\beta} s^2}\frac{ ds}{\sqrt{s^2-k^2}} \\
&=&C \int_k^{\infty} \frac{1}{s^{\alpha+\beta+2} }ds = C(\alpha, \beta) < \infty,
\end{eqnarray}
 therefore
\begin{equation}\label{estimate2}
\int_0^{\pi/2} \frac{sec^{2+\alpha+\beta} \omega \, \cos (\alpha-\beta) \omega}{Z_1^\alpha Z_2^\beta Y} d\omega \leq C(\alpha, \beta) +C k^{-(1+\alpha+\beta)} \int_k^2 \frac{s^{\alpha+\beta+1}}{s^{2\alpha} s^{2\beta} s^2}\frac{ ds}{\sqrt{s^2-k^2}}.
\end{equation}
The Watson kernel is good for localization. The deficits of this representation are:  first, the integral is only convergent for $\alpha + \beta >-1$; second, it is not clear from the representation that the kernel is positive.

There is another representation of the Watson kernel obtained by  W. N. Bailey  in 1939 (\cite{bay}
page 102, see also \cite{askey} page 11),
\begin{eqnarray*}\label{Baileyrep}
K^{(\alpha,\beta)}(r, x, y)&=& \frac{\Gamma( \alpha+\beta+2)(1-r)}{2^{\alpha+\beta+2}\Gamma(\alpha+1)\Gamma(\beta+1)(1+r)^{\alpha+\beta+2}}\\
\nonumber && \quad \quad \quad  \quad \quad \quad  \times \sum_n \sum_m \frac{(\frac{(\alpha+\beta+2)}{2})_{m+n}(\frac{(\alpha+\beta+3)}{2})_{m+n}}{m!n! (\alpha+1)_m(\beta+1)_n} (\frac{a^2}{k^2})^m (\frac{b^2}{k^2})^n\\
&=& \frac{\Gamma( \alpha+\beta+2)(1-r)}{2^{\alpha+\beta+2}\Gamma(\alpha+1)\Gamma(\beta+1)(1+r)^{\alpha+\beta+2}}\\
\nonumber && \quad \quad \quad  \times F_4\left(
\frac{(\alpha+\beta+2)}{2}, \frac{(\alpha+\beta+3)}{2};\alpha+1,\beta+1;\frac{a^2}{k^2},
\frac{b^2}{k^2}\right),\\
\end{eqnarray*}
with $a =\frac{\sqrt{(1-x)(1-y)}}{2}, \, b = \frac{\sqrt{(1+x)(1+y)}}{2},$ and as before $k = \frac{1}{2} (r^{-1/2}+r^{1/2})$. $F_4$ is the Appell
 hypergeometric  function in two variables,
\begin{equation}\label{appellfunc}
 F_4(\alpha,\beta;\gamma, \gamma';x,y) = \sum_n \sum_m \frac{(\alpha)_{m+n}(\beta)_{m+n}}{m!n! (\gamma)_m(\gamma')_n} x^m y^n.
 \end{equation}
 Let us observe that the condition for absolute convergence of the $F_4$ function is $|x|^{1/2} + |y|^{1/2} <1$, see \cite{wangguo}, and therefore the expression above for $K^{(\alpha,\beta)}(r, x, y)$ converges absolutely if $\frac{a}{k}+\frac{b}{k} <1$ and that there is not restriction on $\alpha, \beta$, i.e. it is valid for any $\alpha > -1, \, \beta> -1$.

Moreover by direct inspection of  Bailey's representation it is clear that 
  $$K^{(\alpha,\beta)}(r,x,y) \geq 0.$$ 
 
 From the the uniform convergence of the series of Jacobi polynomials and the fact that the system is complete,  it can be proved, using the orthogonality, that
 \begin{equation}\label{conserv}
\int_{-1}^1 K^{(\alpha,\beta)}(r, x, y) J^{\alpha,\beta}(dy) =1.
\end{equation}
By Holder's inequality, it is easy to see that for $ 1 \leq p \leq \infty$,
\begin{equation}\label{inqp}
\| f^{\alpha, \beta}(r, \cdot)\|_{p,\alpha, \beta} \leq  \| f\|_{p,\alpha, \beta},
\end{equation}
where 
$$  \| f\|_{p,\alpha, \beta} =(\int_{-1}^1 |f(x)|^p J^{\alpha,\beta}(dy))^{1/p},$$
is the $L^p$ norm with respect to the Jacobi measure $ J^{\alpha,\beta}(dy).$ 

Moreover, we have the strong $L^p$-convergence of the Abel sum, we will present an elementary and direct proof of this result.
\begin{lemma}
\begin{equation}
 \| f^{\alpha, \beta}(r, \cdot) -f\|_{p,\alpha, \beta} \rightarrow 0, \quad \mbox{as} \quad r\rightarrow1.
\end{equation}
\end{lemma}

\dem 
\begin{itemize}

\item Using Parserval's identity, the positivity of  $ K^{(\alpha,\beta)}(r, x, y)$ and the completeness of  $\{P^{\alpha,\beta}_n\}$, we have for $f \in L^2(J^{\alpha,\beta})$, 
$$  \| f^{\alpha, \beta}(r, \cdot) -f\|_{2,\alpha, \beta} = \sum_{n=0}^\infty (r^{2n}-1) | \hat{f}^{(\alpha,\beta)}(n)|^2 \rightarrow 0,$$
as $r \rightarrow 1$.

For the other cases $p \neq 2$ given $\lambda >0$ fix, and $f \in L^p(J^{\alpha,\beta}),$ without loss of generality we may assume $f \geq 0$ and then we can write $f$  as $f = f_1 + f_2$ with $|f_1| \leq \lambda$, $ f_1 \in  L^2(J^{\alpha,\beta})$ and let us take $\lambda$ big  enough that $\|f_2\|_p< \varepsilon$.

\item Now if $2 < p \leq \infty$, then 
 $|\frac{f_1}{\lambda}| \leq 1$ implies $ |\frac{f_1}{\lambda}|^p \leq |\frac{f_1}{\lambda}|^2$,

\begin{eqnarray*}
 \| f_1^{\alpha, \beta}(r, \cdot) -f_1\|^p_{p,\alpha, \beta}& = & 2^p \lambda^p\|\frac{1}{2} (\frac{f_1}{\lambda})^{\alpha, \beta}(r, \cdot) - \frac{1}{2}(\frac{f_1}{\lambda})\|^p_{p,\alpha, \beta}\\
 & \leq & 2^p \lambda^p\|\frac{1}{2} (\frac{f_1}{\lambda})^{\alpha, \beta}(r, \cdot) -\frac{1}{2}(\frac{f_1}{\lambda})\|^2_{2,\alpha, \beta} \\
 & = & 2^{p-2} \lambda^{p-2}\|f_1^{\alpha, \beta}(r, \cdot) -f_1\|^2_{2,\alpha, \beta} \rightarrow 0\end{eqnarray*}
as $r \rightarrow 1$, from the previous case. Now from (\ref{inqp})
$$ \| f_2^{\alpha, \beta}(r, \cdot) -f_2\|^p_{p,\alpha, \beta} \leq  2^p (\| f_2^{\alpha, \beta}(r, \cdot)\|^p_{p,\alpha, \beta} +\|f_2\|^p_{p,\alpha, \beta}) \leq  2^{p+1}\|f_2\|^p_{p,\alpha, \beta} < 2^{p+1} \varepsilon^p.$$

\item Finally, for $1\leq p <2$, from  (\ref{conserv}) (taking $s>1$ such that $ sp =2$) and using H\"older's inequality,
\begin{eqnarray*}
 \| f_1^{\alpha, \beta}(r, \cdot) -f_1\|^p_{p,\alpha, \beta}&\leq& C  \| f_1^{\alpha, \beta}(r, \cdot) -f_1\|^2_{2,\alpha, \beta}. 
\end{eqnarray*}
The inequality for $f_2$ is obtained similarly as in the previous case.
\end{itemize}

\ep

The Jacobi maximal function $f^*_{\alpha, \beta}$, is defined as
\begin{equation}\label{JacobMaxfunt}
f^*_{\alpha, \beta}(x) =  \sup_{0<r<1} |f^{\alpha, \beta}(r,x)| =  \sup_{0<r<1} | \int_{-1}^{1} K^{\alpha,\beta}(r,x,y) f(y) J^{\alpha,\beta}(dy)|. 
\end{equation}
We will prove, as a consequence of the main result of this paper, that $f^*_{\alpha, \beta}$ is weak-$(1,1)$  continuous with respect to $J^{\alpha,\beta}$, i. e.
\begin{eqnarray}\label{weak1}
J^{\alpha, \beta} \{ f_{\alpha, \beta}^* >\lambda\} \leq \frac{C_{\alpha, \beta}}{\lambda} \| f\|_{1,\alpha, \beta}.
\end{eqnarray}
From Bailey's representation it is almost trivial to get
\begin{equation}
\| f^{\alpha, \beta}(r, \cdot)\|_{\infty} \leq C \| f\|_{\infty},
\end{equation}
then
\begin{equation}
\| f^*_{\alpha, \beta}(r, \cdot)\|_{\infty} \leq C \| f\|_{\infty},
\end{equation}
therefore, by interpolation we get, for $1 < p < \infty,$
\begin{equation}
\| f^*_{\alpha, \beta}(r, \cdot)\|_{p,\alpha, \beta} \leq C \| f\|_{p,\alpha, \beta}
\end{equation}

For more details on the Jacobi maximal function can be found in \cite{cafacal}, \cite{cafacal2} and \cite{calv}.

\section{Estimates of the Watson kernel}
By the product rule in the Watson representation (\ref{Watsonrep}),
$$K^{\alpha,\beta}(r,x,y) = r^{(1-\alpha-\beta)/2} \frac{d}{dr}\left( k^{1+\alpha+\beta} \int_0^{\pi/2} \frac{sec^{2+\alpha+\beta} \omega \, \cos (\alpha-\beta) \omega}{Z_1^\alpha Z_2^\beta Y} d\omega\right)$$
 we get four kernels $A, B, C, D$
defined in the following way,
\begin{eqnarray*}
A&=&  r^{(1-\alpha-\beta)/2} \frac{d}{dr}(k^{1+\alpha+\beta}) \int_0^{\pi/2} \frac{sec^{2+\alpha+\beta} \omega \, \cos (\alpha-\beta) \omega}{Z_1^\alpha Z_2^\beta Y} d\omega,\\
B&=&  r^{(1-\alpha-\beta)/2} k^{1+\alpha+\beta} \int_0^{\pi/2} \frac{d}{dr}(Y^{-1})\frac{sec^{2+\alpha+\beta} \omega \, \cos (\alpha-\beta) \omega}{Z_1^\alpha Z_2^\beta} d\omega,\\
C &=&  r^{(1-\alpha-\beta)/2} k^{1+\alpha+\beta} \int_0^{\pi/2} \frac{d}{dr}(Z_1^{-\alpha})\frac{sec^{2+\alpha+\beta} \omega \, \cos (\alpha-\beta) \omega}{ Z_2^\beta Y} d\omega,\\
D&=&  r^{(1-\alpha-\beta)/2} k^{1+\alpha+\beta} \int_0^{\pi/2} \frac{d}{dr}(Z_2^{-\beta})\frac{sec^{2+\alpha+\beta} \omega \, \cos (\alpha-\beta) \omega}{ Z_1^{\alpha} Y} d\omega.
\end{eqnarray*}
Then we have, see \cite{cafacal2} pages 282-3 or \cite{calv} Lemma 4.1, pages 245-9,
\begin{lemma}
We have the following estimate for the Watson kernel,
\begin{equation}\label{basicineq}
K^{\alpha,\beta}(r,x,y)  \leq C(\alpha, \beta) (1+L(r,x,y)),
\end{equation}
where $C(\alpha,\beta)$ is a positive constant, $L(r,x,y)$ is the integral
\begin{equation}\label{likepoisson}
L(r,x,y)=(1-r) \int_k^{2} \frac{(s-min(x,y))^{1-\alpha}}{((x-y)^2+(s-1)(s-min(x,y)))^{3/2}} \frac{ds}{(s-k)^{1/2}},
\end{equation}
where $ k = \frac{1}{2} (r^{1/2} + r^{-1/2}), \,0 \leq x \leq 1.$
\end{lemma}

For the proof of this lemma, the following estimates  will be needed, for detail see Appendix in \cite{calv}. Let $1\leq s \leq 2, 0 \leq x \leq 1, |y| \leq 1.$ Then:
 \begin{enumerate}
\item [i)] $ s^2 -\min(x,y) \leq 4(s- \min(x,y));$
\item [ii)] $ s -\min(x,y) \leq 2(s-xy) \leq 4(s- \min(x,y));$
\item[iii)] $C_1\left((x-y)^2+ (s-1)(s-\min(x,y)) \right) \leq Y^2 \leq C_2 \left((x-y)^2+ (s-1)(s-\min(x,y)) \right);$
\item [iv)] $ s^2 -\min(x,y) \leq Z_1 \leq C (s^2 -\min(x,y));$
\item[v)] $1 \leq s^2 + \max(x,y) \leq Z_2 \leq C;$
\item[vi)] If $\varphi(x,r) = (k-1)^{1/2}(k-x)^{1/2}$, then $k-1 \leq \phi(x,r) \leq k-x,$ for $k > 1$;
\item[vii)] $C_1 (1-r)^2 \leq k-1 \leq C_2 (1-r^2)$, if $0 < r_0 < r < 1.$
\end{enumerate}
Here $C, C_1, C_2$ denote positive constants.
 From these estimates observe that:
 \begin{itemize}
\item By iii), $Y^2 \sim  \left((x-y)^2+ (s-1)(s-\min(x,y)) \right).$
\item By iv),  $Z_1 \sim  (s^2 -\min(x,y)).$
\item By v), $Z_2$ is essentially a constant.
\end{itemize}

Observe that if $-1 < x <0$ similar estimates hold, just changing the role of $\alpha$ and $\beta$.
For details of the proof of Lemma 1 see \cite{calv} Lemma 4.1.

 In  \cite{cafacal2} pages 284-6 and \cite{calv} Lemma 4.1, page 254, the following estimate for $L$ 
was obtained,
\begin{lemma}
\begin{equation}
L(r,x,y) \leq C_{\alpha,\beta} \sum_{n=0}^{\infty} \frac{1}{2^{n/2}} \frac{1}{J^{\alpha,\beta}(I_n(x,r))}\chi_{I_n(x,r)},
\end{equation}
where $I_n(x,t) = [x - 2^n \varphi(x,r), x + 2^n \varphi(x,r)] \cap [-1,1]$, $\chi_{I_n(x,r)}$ is its characteristic function and $\varphi(x,r) = (k-1)^{1/2} (k-x)^{1/2}.$
\end{lemma}

We are going to get another estimate related to $L(r,x,y)$ using superposition of  Poisson  type kernels.
The following technical result, see $(5.1)$ and $(5.2)$ of \cite{calv}, is needed,
\begin{lemma} There exist constants $C_1$ and $C_2$ independent of r< such that,
 \begin{equation}\label{estm1}					
(1-r) \int_k^2 \frac{1}{(s-k)^{1/2}(s-x)^{1/2}} ds < C_1
\end{equation}
and 
 \begin{equation}\label{estm2}
(1-r) \int_k^2 \frac{1}{(s-k)^{1/2}( s-1)^{1/2}(s-x)^{1/2}} ds < C_1
\end{equation}
\end{lemma}
\dem

Let us prove first (\ref{estm1}). Observe that, by the estimate vii) we have $(k-1) \sim (1-r)^2$ i.e. $(k-1)^{1/2} \sim (1-r)$. Then, integrating by parts,
\begin{eqnarray*}
(k-1)^{1/2} \int_k^2 \frac{1}{(s-k)^{1/2}(s-1)} ds &=& (k-1)^{1/2} [\frac{2(s-k)^{1/2}}{(s-1)} ]_k^2 + 2 \int_k^2 \frac{(s-k)^{1/2}}{(s-1)^2}ds]\\
&=& (k-1)^{1/2} [2(2-k)^{1/2} + \int_k^2 \frac{(s-k)^{1/2}}{(s-1)^2}ds],\\
\end{eqnarray*}
and 
\begin{eqnarray*}
(k-1)^{1/2}\int_k^2 \frac{(s-k)^{1/2}}{(s-1)^2}ds &\leq&(k-1)^{1/2} \int_k^2 \frac{1}{(s-1)^{3/2}}ds\\
&=&(k-1)^{1/2}\int_k^2 \frac{1}{(s-k+k-1)^{3/2}}ds\\
&\leq& \int_k^2 \frac{1}{(k-1)}\frac{1}{(|\frac{s-k}{k-1}|+1)^{3/2}}ds\\
&=& \frac{1}{\lambda} \int_k^2 k_1(\frac{s-k}{\lambda}) ds <C,
\end{eqnarray*}
it $\lambda =(k-1)$ and the the Poisson type kernel $k_1(x) = \frac{1}{ (|x|+1)^{3/2}}$. Observe that 
$\int_{-\infty}^{\infty} k_1(x) dx = \int_{-\infty}^{\infty} \frac{1}{ (|x|+1)^{3/2} }dx=4.$\\

The second estimate (\ref{estm2}) follows immediately from (\ref{estm1}).
\ep\\

The following technical result is also needed for the proof of Theorem \ref{mainest}, 

\begin{lemma}\label{kernelest} For any $\eta>1 $
$$ \sup_{0<|a| <1} \frac{1}{[(z+a)^2+1]^{\eta}} \leq  \frac{C}{[z^2+1]^{\eta}}.$$
\end{lemma}
\dem

-  If $|z| >3$ i.e. $\frac{|z|}{3} >1$, then for $0 < |a| <1$
\begin{eqnarray*}
|z+a| &\geq& |z| - |a| > \frac{2|z|}{3}+ (\frac{|z|}{3} -1) \geq  \frac{2|z|}{3}, \,  \text{so}\\
|z+a|^2 &\geq& \frac{4|z|^2}{9}.
\end{eqnarray*}
Thus 
\begin{eqnarray*}
 \frac{1}{[(z+a)^2+1]^{\eta}} &\leq& \frac{1}{[ \frac{4|z|^2}{9}+1]^{\eta}}\\
 & \leq&  \frac{1}{[ \frac{4|z|^2}{9}+\frac{4}{9}]^{\eta}} =\frac{(\frac{9}{4})^{\eta}}{[ |z|^2+1]^{\eta}} = \frac{C}{[ |z|^2+1]^{\eta}}.
 \end{eqnarray*}
 - If $|z|< 3$, then
$$ \frac{1}{[(z+a)^2+1]^{\eta}} \leq 1, \quad  
\text{and} \quad \frac{1}{10^\eta} \leq \frac{1}{[z^2+1]^{\eta}} \leq 1,$$
thus
$$ \frac{1}{[(z+a)^2+1]^{\eta}} \leq 1 \leq  \frac{10^\eta}{[z^2+1]^{\eta}} = \frac{C}{[z^2+1]^{\eta}} .$$
 \ep

The main estimates of the Watson kernel that we have obtained in this paper is the following,

\begin{theorem}\label{mainest}
The integral,
\begin{equation}\label{intest}
\int_0^1(1-r) \int_k^{2} \frac{(s-\min(x,y))^{1-\alpha}}{((x-y)^2+(s-1)(s- \min(x,y)))^{3/2}} \frac{ds}{(s-k)^{1/2}} (1-y)^{\alpha} \, dy
\end{equation}
is bounded by a superposition of a family of Poisson type kernels integrated with respect to a parameter, and therefore it is bounded from above.
\end{theorem}

\dem
\begin{enumerate}
\item[i)]  Case $\alpha\geq 0.$
\begin{enumerate}
\item[ i-1)] If $x \leq y <1$: in this range $(s-x)^{-\alpha} \leq (1-y)^{-\alpha}, $ and hence
\begin{eqnarray*}
&&\int_x^1(1-r) \int_k^{2} \frac{(s-x)(s-x)^{-\alpha}}{((x-y)^2+(s-1)(s-x))^{3/2}} \frac{ds}{(s-k)^{1/2}} (1-y)^{\alpha} \, dy\\
&\leq& (1-r) \int_k^{2} \frac{1}{(s-k)^{1/2}} \frac{(s-x)}{(s-x)(s-1)}\\
&&\quad \quad \quad \quad \times \int_x^1  \frac{1}{[(s-1)(s-x)]^{1/2}} \frac{dy}{((\frac{x-y}{[(s-1)(s-x)]^{1/2}})^2+1)^{3/2}} ds
\end{eqnarray*}
Considering the Poisson type kernel  $k_2(x) = \frac{1}{(x^2+1)^{3/2}}$ then the inner integral can be rewritten as
$$  \int_x^1  \frac{1}{[(s-1)(s-x)]^{1/2}} \frac{dy}{((\frac{x-y}{[(s-1)(s-x)]^{1/2}})^2+1)^{3/2}} ds
=\frac{1}{\lambda} \int_x^1k_2(\frac{x-y}{\lambda}) dy,$$
with $\lambda= [(s-1)(s-x)]^{1/2},$ and since $\int_{-\infty}^{\infty} k_2(x) dx = \int_{-\infty}^{\infty}  \frac{1}{(x^2+1)^{3/2}} dx =2,$ then the  inner integral is bounded and therefore
\begin{eqnarray*}
&&\int_x^1(1-r) \int_k^{2} \frac{(s-x)(s-x)^{-\alpha}}{((x-y)^2+(s-1)(s-x))^{3/2}} \frac{ds}{(s-k)^{1/2}} (1-y)^{\alpha} \, dy\\
&\leq& C (1-r) \int_k^2 \frac{1}{(s-k)^{1/2}(s-1)} ds < C,
\end{eqnarray*}
by $(5.2)$ of \cite{calv}.
\item[ i-2)] If $0< y< x$: in this range 
$$(s-\min(x,y))^{-\alpha} (1-y)^{\alpha}=(s-y)^{-\alpha} (1-y)^{\alpha} \leq 1,$$ 
  then the corresponding part of (\ref{intest}) in this range is less than
\begin{eqnarray*}
&&\int_0^x(1-r) \int_k^{2} \frac{(s-y)}{((x-y)^2+(s-1)(s-y))^{3/2}} \frac{ds}{(s-k)^{1/2}}  \, dy\\
&\leq & \int_0^x(1-r) \int_k^{2} \frac{(s-y)}{((x-y)^2+(s-1)(s-x))^{3/2}} \frac{ds}{(s-k)^{1/2}}  \, dy\\
&\leq & (1-r) \int_k^{2} \frac{1}{(s-k)^{1/2}} \frac{(s-y)}{(s-1)(s-x)}\\
&&\quad \quad \quad \quad \times \int_0^x  \frac{1}{(s-1)^{1/2}(s-x)^{1/2}} \frac{dy}{((\frac{x-y}{[(s-1)(s-x)]^{1/2}})^2+1)^{3/2}} ds
\end{eqnarray*}
Now as $s-y = s-x+x-y,$ we get two terms,
\begin{eqnarray*}
& & (1-r) \int_k^{2} \frac{1}{(s-k)^{1/2}} \frac{1}{(s-1)}\\
&&\quad \quad \quad \quad \times \int_0^x  \frac{1}{(s-1)^{1/2}(s-x)^{1/2}} \frac{dy}{((\frac{x-y}{[(s-1)(s-x)]^{1/2}})^2+1)^{3/2}} d\\
&+& (1-r) \int_k^{2} \frac{1}{(s-k)^{1/2}} \frac{1}{(s-1)^{1/2}(s-x)^{1/2}}\\
&&\quad \quad \quad \quad \times \int_0^x  \frac{1}{(s-1)^{1/2}(s-x)^{1/2}} \frac{[\frac{x-y}{(s-1)^{1/2}(s-x)^{1/2}}]}{((\frac{x-y}{[(s-1)(s-x)]^{1/2}})^2+1)^{3/2}} \,dy\,ds
\end{eqnarray*}
The first integral is analogous to case i-1) i.e it is bounded by
$$ C (1-r) \int_k^2 \frac{1}{(s-k)^{1/2}(s-1)} ds < C_1,$$
by (\ref{estm1}). The second integral is bounded by
\begin{eqnarray*}
&& (1-r) \int_k^{2} \frac{1}{(s-k)^{1/2}} \frac{1}{(s-1)^{1/2}(s-x)^{1/2}}\\
&&\quad \quad \quad \quad \times \int_0^x  \frac{1}{(s-1)^{1/2}(s-x)^{1/2}} \frac{([\frac{x-y}{(s-1)^{1/2}(s-x)^{1/2}}]^2+1)^{1/2}}{((\frac{x-y}{[(s-1)(s-x)]^{1/2}})^2+1)^{3/2}} \,dy\,ds\\
&=& (1-r) \int_k^{2} \frac{1}{(s-k)^{1/2}} \frac{1}{(s-1)^{1/2}(s-x)^{1/2}}\\
&&\quad \quad \quad \quad \times \int_0^x  \frac{1}{(s-1)^{1/2}(s-x)^{1/2}} \frac{1}{((\frac{x-y}{[(s-1)(s-x)]^{1/2}})^2+1)} \,dy\,ds,
\end{eqnarray*}
and therefore we get the bound
$$ C (1-r) \int_k^2 \frac{1}{(s-k)^{1/2}(s-k)^{1/2}(s-1)^{1/2}} ds < C_2,$$
by considering the Poisson type kernel $ k_3  = \frac{1}{(x^2+1)}$, as $\int_{-\infty}^{\infty} k_3(x) dx = \int_{-\infty}^{\infty}  \frac{1}{(x^2+1)} dx =\pi,$ and estimate (\ref{estm2}).\\
\end{enumerate}
\item[ii)]  Case $-1<\alpha <0.$
\begin{enumerate}
\item[ ii-1)] If $x \leq y <1$: we rewrite the corresponding part of (\ref{intest}) in this range as
\begin{eqnarray*}
&&(1-r) \int_k^{2}\frac{1}{(s-k)^{1/2}(s-1)(s-x)} \\
&&\quad \quad \quad \quad \times \int_x^1  \frac{(s-x)^{1-\alpha}(1-y)^{\alpha} }{(s-1)^{1/2}(s-x)^{1/2}} \frac{dy}{((\frac{x-y}{[(s-1)(s-x)]^{1/2}})^2+1)^{3/2}}  \, ds\\
\end{eqnarray*}
Now
$$ x-y = [(x+1 -s)  -y] +(s -1),$$
hence
\begin{eqnarray*}
\frac{x-y}{[(s-1)(s-x)]^{1/2}} = \frac{(x+1 -s)  -y}{[(s-1)(s-x)]^{1/2}} +\frac{s-1}{[(s-1)(s-x)]^{1/2}},
\end{eqnarray*}
 then, if $a = \frac{s-1}{[(s-1)(s-x)]^{1/2}}$, the inner integral can be rewritten as
\begin{eqnarray*}
&&(s-x)^{-\alpha} \int_x^1 \frac{(1-y)^{\alpha} }{(s-1)^{1/2}(s-x)^{1/2}} \frac{dy}{[(\frac{(x+1-s)-y}{[(s-1)(s-x)]^{1/2}}+a)^2+1]^{3/2}},\\
\end{eqnarray*}
then by the Lemma \ref{kernelest},  with $\eta= 3/2$, we get
$$ \sup_{0<a <1} \frac{1}{[(z+a)^2+1]^{3/2}} \leq  \frac{C}{[z^2+1]^{3/2}}.$$
Therefore, we get
\begin{eqnarray*}
&&C (1-r) \int_k^{2}\frac{(s-x)^{-\alpha}}{(s-k)^{1/2}(s-1)} \\
&&\quad \quad \quad \quad \times \int_x^1  \frac{(1-y)^{\alpha} }{(s-1)^{1/2}(s-x)^{1/2}} \frac{dy}{[(\frac{(x+1 -s)  -y}{[(s-1)(s-x)]^{1/2}})^2+1]^{3/2}}  \, ds\\
\end{eqnarray*}

Considering again the  Poisson type kernel $k_2(x) = \frac{1}{(x^2+1)^{3/2}}$ and $\lambda = [(s-1)(s-x)]^{1/2}$ this can be written as
\begin{eqnarray*}
&&C (1-r) \int_k^{2}\frac{(s-x)^{-\alpha}}{(s-k)^{1/2}(s-1)} \\
&&\quad \quad \quad \quad \times \frac{1}{\lambda} \int_x^1  (1-y)^{\alpha} k_2( \frac{(x+1 -s)  -y}{\lambda})dy  \, ds.\\
\end{eqnarray*}
By a classical argument the inner integral in previous expression is bounded by  $M \psi(x+1-s)$ where $M \psi$ is the Hardy-Littlewood maximal function of  $\psi(y) = (1-y)^{\alpha}$. Now since $\psi$ is a $A_1$- Muckenhoupt weight with respect to the Lebesgue measure, see \cite{duo}, we get that the inner integral is then bounded by
$$  M \psi(x+1-s) \leq C  \psi(x+1-s)= C [1-(x+1-s)]^{\alpha} = C  (s-x)^{\alpha} .$$
Thus,
the corresponding part of (\ref{intest}) in this range is bounded by
$$ C (1-r) \int_k^2 \frac{1}{(s-k)^{1/2}(s-1)} ds < C,$$
by estimate (\ref{estm1}).
\item[ ii-2)] If $0< y< x$: The corresponding part of (\ref{intest}) in this range takes de form,
\begin{eqnarray*}
&&\int_0^x(1-r) \int_k^{2} \frac{(s-y)^{1-\alpha}}{((x-y)^2+(s-1)(s-y))^{3/2}} \frac{ds}{(s-k)^{1/2}} (1-y)^{\alpha} \, dy\\
&\leq & (1-r) \int_k^{2} \frac{1}{(s-k)^{1/2}(s-1)(s-x)}\\
&&\quad \quad \quad \quad \times \int_0^x  \frac{(s-y)^{1-\alpha}(1-y)^{\alpha}}{(s-1)^{1/2}(s-x)^{1/2}((\frac{x-y}{[(s-1)(s-x)]^{1/2}})^2+1)^{3/2}} dy \, ds.
\end{eqnarray*}
Now since $\alpha <0$
$$(s-y)^{1-\alpha} \leq C_\alpha [(s-x)^{1-\alpha}+(x-y)^{(1-\alpha)}],$$
 we get two terms
\begin{eqnarray*}
&&C_\alpha  (1-r) \int_k^{2} \frac{(s-x)^{1-\alpha}}{(s-k)^{1/2}(s-1)(s-x)}\\
&&\quad \quad \quad \quad \times \int_0^x  \frac{(1-y)^{\alpha}}{(s-1)^{1/2}(s-x)^{1/2}((\frac{x-y}{[(s-1)(s-x)]^{1/2}})^2+1)^{3/2}} dy \, ds \\
& &\quad \quad \quad + C_\alpha (1-r) \int_k^{2} \frac{1}{(s-k)^{1/2}(s-1)(s-x)}\\
&&\quad \quad \quad \quad \times \int_0^x  \frac{(x-y)^{1-\alpha}(1-y)^{\alpha}}{(s-1)^{1/2}(s-x)^{1/2}((\frac{x-y}{[(s-1)(s-x)]^{1/2}})^2+1)^{3/2}} dy \, ds
\end{eqnarray*}
The first integral can be handle in a similar way as in the case ii-1); taking
$$ x-y = [(x+1 -s)  -y] +(s -1),$$
and using again Lemma \ref{kernelest},  with $\eta= 3/2$, we get as before,
\begin{eqnarray*}
&& (1-r) \int_k^{2} \frac{(s-x)^{-\alpha}}{(s-k)^{1/2}(s-1)}\\
&&\quad \quad \quad \quad \times \int_0^x  \frac{(1-y)^{\alpha}}{(s-1)^{1/2}(s-x)^{1/2}((\frac{x-y}{[(s-1)(s-x)]^{1/2}})^2+1)^{3/2}} dy \, ds.
\end{eqnarray*}
 Then the inner integral is less or equal than $C(s-x)^\alpha$ and therefore this term is less than
$$ C(1-r) \int_k^{2} \frac{1 }{(s-k)^{1/2}(s-1)} ds < C,$$
using estimate (\ref{estm1}).

For the second integral the numerator of the inner integral can be rewritten as
\begin{eqnarray*}
(x-y)^{1-\alpha} &=& [(s-1)^{1/2} (s-x)^{1/2}]^{1-\alpha}(\frac{x-y}{(s-1)^{1/2}(s-x)^{1/2})})^{1-\alpha}\\
&=& [(s-1)^{1/2} (s-x)^{1/2}]^{1-\alpha}[(\frac{x-y}{(s-1)^{1/2}(s-x)^{1/2})})^2]^{(1-\alpha)/2}\\
&\leq& [(s-1)^{1/2} (s-x)^{1/2}]^{1-\alpha}[(\frac{x-y}{(s-1)^{1/2}(s-x)^{1/2}})^2+1]^{(1-\alpha)/2}\\
\end{eqnarray*}

Then the inner integral  is bounded by
$$ \int_0^x  \frac{1}{(s-1)^{1/2}(s-x)^{1/2}}\frac{1}{[(\frac{x-y}{[(s-1)(s-x)]^{1/2}})^2+1]^{3/2-(1-\alpha)/2}} (1-y)^{\alpha}dy,$$
and $\frac{3}{2} - \frac{1-\alpha}{2}= 1+\alpha/2 >1/2.$
Therefore the second integral is bounded by
\begin{eqnarray*}
&& (1-r) \int_k^{2} \frac{(s-1)^{(1-\alpha)/2} (s-x)^{(1-\alpha)/2}}{(s-k)^{1/2}(s-1)(s-x)}\\
&&\quad \quad \times \int_0^x  \frac{1}{(s-1)^{1/2}(s-x)^{1/2}}\frac{1}{[(\frac{x-y}{[(s-1)(s-x)]^{1/2}})^2+1]^{3/2-(1-\alpha)/2}} (1-y)^{\alpha}dy \, ds\\
&& \leq  (1-r) \int_k^{2} \frac{(s-x)^{(1-\alpha)/2} (s-x)^{(1-\alpha)/2}}{(s-k)^{1/2}(s-1)(s-x)}\\
&&\quad \quad \times \int_0^x  \frac{1}{(s-1)^{1/2}(s-x)^{1/2}}k_4(\frac{x-y}{(s-1)^{1/2}(s-x)^{1/2}}) (1-y)^{\alpha}dy \, ds\\
&& \leq  (1-r) \int_k^{2} \frac{(s-x)^{-\alpha}}{(s-k)^{1/2}(s-1)}\\
&&\quad \quad \times \int_0^x  \frac{1}{(s-1)^{1/2}(s-x)^{1/2}}k_4(\frac{x-y}{(s-1)^{1/2}(s-x)^{1/2}}) (1-y)^{\alpha}dy \, ds.
\end{eqnarray*}

Then, this is analogous to the case ii-1), but with the Poisson type kernel $k_4(x) = \frac{1}{(x^2+1)^{3/2-(1-\alpha)/2}},$ and $\lambda =(s-1)^{1/2}(s-x)^{1/2}.$\\
and therefore the second integral is bounded by
$$(1-r) \int_k^2 \frac{1}{(s-k)^{1/2}(s-1)} ds < C,$$
using estimate (\ref{estm1}).
\end{enumerate}

\ep
\end{enumerate}

\section{Applications}
We are going to obtain several consequences from Theorem \ref{mainest}.\\

First we need to consider a result due to A. Zygmund (see \cite{zy} Vol I Lemma 7.1 page 154-5) which in particular implies Natanson's  lemma (see  \cite{mu1} Theorem 1),
\begin{lemma}(Zygmund) \label{zyglem} 

Given $-\infty \leq a < b \leq \infty$ a Borel measure $\mu$ with support in $(a,b)$  and a kernel  $K(r, x,\cdot)$ depending of a parameter $r$, such that 
\begin{equation}\label{Zygbound1}
 \int_a^b |K(r,x,y)| \mu(dy) \leq M_1
\end{equation}
and 
\begin{equation}\label{Zygbound2}
 \int_x^b \mu(x,y) V_2(K(r, x,dy)) \leq M_2, \quad \int_a^x  \mu(y,x) V_2(K(r, x,dy)) \leq M_2,
\end{equation}

where $M_1, M_2$ are constants independent of $x$ and $r$, $V_2(K(r, x,\cdot))$ is the (first) variation of the kernel $K(r,x,y) $ in the variable $y$, i. e. 
$$ V_2(K(r, x,\cdot)) = \sup \sum_i |K(r, x,y_i) - K(r, x,y_{i-1})|,$$
where the supremum is taken over all partitions of $[a,b]$ and the integrals are considered in the Lebesgue-Stieltjes sense.

Then for $f \in L^1(\mu) $,
\begin{equation}
|\int_a^b K(r, x,y) f(y) \mu(dy)| \leq M f_\mu^*(x),
\end{equation}
where $M$ depends only on $M_1, M_2$ and
$$ f_\mu^*(x) = \sup_{x \in I} \frac{1}{\mu(I)} \int_I f(y) \mu(dy),$$
 is the non-centered Hardy-Littlewood  maximal  function for $f$ with respect to the measure  $\mu$. 
 \end{lemma}
 
 \dem

Using the integration by parts formula for Stieltjes integrals, we have

\begin{eqnarray*}
 \int_x^b K(r,x,y) \mu(dy)& =& ( \int_x^b \mu(du) ) K(r,x,b) - \int_x^b(\int_x^y \mu(du) ) K(r,x,dy)\\
 & =& \mu(x,b) K(r,x,b)- \int_x^b \mu(x,y)  K(r,x,dy).
\end{eqnarray*}
Therefore, by hypothesis
\begin{eqnarray*}
  |\mu(x,b) K(r,x,b) | &\leq&  \int_x^b |K(r,x,y)| \mu(dy) +  \int_x^b \mu(x,y)  K(r,x,dy)\\
  &\leq&  \int_x^b |K(r,x,y)| \mu(dy) +  \int_x^b \mu(x,y)  V_2(K(r,x,dy))\\
   &\leq& M_1 +M_2
\end{eqnarray*}
Now, for $f \in L^1(\mu) $ using again the integration by parts formula,
\begin{eqnarray*}
 \int_x^b f(y) K(r,x,y) \mu(dy)& =& ( \int_x^b f(y) \mu(dy)) K(r,x,b) - \int_x^b(\int_x^y f(y) \mu(dy) ) K(r,x,dy)\\
 & =& ( \int_x^b f(y) \mu(dy)) K(r,x,b)- \int_x^b( \int_x^b f(y) \mu(dy)) K(r,x,dy)\\
 & =& ( \frac{1}{\mu(x,b)} \int_x^b f(y) \mu(dy))\mu(x,b) K(r,x,b) \\
 && \quad \quad  \quad  \quad \quad \quad  - \int_x^b( \frac{1}{\mu(x,y)}  \int_x^b f(y) \mu(dy))\mu(x,y) K(r,x,dy).
\end{eqnarray*}
Thus,
\begin{eqnarray*}
| \int_x^b f(y) K(r,x,y) \mu(dy)|& \leq & f_\mu^*(x) |\mu(x,b) K(r,x,b)| +  f_\mu^*(x)  \int_x^b\mu(x,y) V_2(K(r,x,dy))\\
&\leq & (M_1+M_2) f_\mu^*(x) + M_2 f_\mu^*(x) = (M_1+ 2 M_2) f_\mu^*(x). 
\end{eqnarray*}
\ep  

\begin{obs}
Given a measure $\mu$ as before, observe that for a Natanson's kernel $K(r,x,y)$( i.e.$-\infty \leq a < b \leq \infty$  and  $K(r, x,y)$ non-negative, such that $K(r,x,y)$ is monotone increasing for $a<y < x $ and monotone decreasing for $b> y >x$, and  $ \int_x^b K(r, x,y) dy =M_1, \quad \int_a^x K(r, x,y) dy =M_2,$
where $M_1, M_2$ are constants independent of $x$), then  $K$ satisfies the conditions of Zygmund's lemma since (\ref{Zygbound1}) is trivial and (\ref{Zygbound2}) is easy obtained by monotonicity conditions. In particular, Poisson type kernels satisfy  the conditions of Zygmund's lemma.\\
\end{obs}

Now, as a  consequence of Theorem \ref{mainest} and using Zygmund's lemma we have,

\begin{theorem}
Let $f \in L^1(J^{\alpha,\beta})$, then the operator
\begin{eqnarray}
\nonumber J_\alpha f(x) &=&\int_0^1(1-r) \int_k^{2} \frac{(s-\min(x,y))^{1-\alpha}}{((x-y)^2+(s-1)(s- \min(x,y)))^{3/2}} \frac{ds}{(s-k)^{1/2}} \\
&& \quad \quad \quad  \quad \quad \quad \quad \quad  \quad \quad \quad  \quad \quad \quad \quad  \quad \quad \quad  \times(1-y)^{\alpha} f(y) \, dy.
\end{eqnarray}
Then,
 \begin{equation}
J_\alpha f(x) \leq  C f^*_{J^{\alpha,\beta}}(x),
\end{equation}
where $f^*_{J^{\alpha,\beta}}$ is the (non-centered) Hardy-Littlewood maximal function with respect to the Jacobi measure
$J^{\alpha,\beta}$. 
 \end{theorem}
\dem

The idea of the proof is the following: by Theorem \ref{mainest} if $f\equiv 1$ for the case $-1 < \alpha <0$ as well as for the case $0 \leq \alpha$ we know that $J_\alpha f(x)$ is bounded by Poisson type kernels and therefore bounded, then using Zygmund's lemma  for the Poisson type kernels we get the result with the (non-centered) Hardy-Littlewood maximal function. The Poisson type kernels are the same used in the proof of Theorem \ref{mainest}

We need to analyze two cases:
\begin{enumerate}
\item[i)]  Case $\alpha\geq 0.$
\begin{enumerate}
\item[ i-1)] If $y> x$: in this range we have,
\begin{eqnarray*}
&&J_\alpha f(x)  \leq (1-r) \int_k^{2} \frac{1}{(s-k)^{1/2}(s-1)}  \int_x^1  \frac{1}{(s-x)^\alpha}\\
&&\quad \quad \quad \quad  \quad \quad  \quad  \times  \frac{1}{[(s-1)(s-x)]^{1/2}} \frac{1}{((\frac{x-y}{[(s-1)(s-x)]^{1/2}})^2+1)^{3/2}} f(y) (1-y)^\alpha dy \, ds
\end{eqnarray*}
Then, 
$$  \frac{1}{[(s-1)(s-x)]^{1/2}} \frac{1}{((\frac{x-y}{[(s-1)(s-x)]^{1/2}})^2+1)^{3/2}}  =  \frac{1}{\lambda} k_3(\frac{x-y}{\lambda}) $$
been a Poisson type kernel, the expression in the inner integral satisfies a (unilateral) condition of  Zygmund's lemma with respect to the measure $\mu(dy) = (1-y)^\alpha dy$ and therefore
\begin{eqnarray*}
&&J_\alpha f(x)  \leq C f^*_{J^{\alpha,\beta}}(x) (1-r) \int_k^{2} \frac{1}{(s-k)^{1/2}(s-1)}  \int_x^1  \frac{1}{(s-x)^\alpha} \frac{1}{\lambda} k_2(\frac{x-y}{\lambda}) (1-y)^\alpha dy \, ds,
\end{eqnarray*}
where $\lambda= [(s-1)(s-x)]^{1/2}$. Then by  the proof of  Theorem \ref{mainest} i-1) we get the last term is bounded i.e.
$$ J_\alpha f(x)  \leq C f^*_{J^{\alpha,\beta}}(x).$$
\item[ i-2)] If $y\leq x$: in this range we have,
\begin{eqnarray*}
&&J_\alpha f(x)  \leq (1-r) \int_k^{2} \frac{1}{(s-k)^{1/2}(s-1)(s-x)}  \int_0^x (s-y)^{1-\alpha}\\
&&\quad \quad \quad \quad  \quad \quad  \quad  \times  \frac{1}{[(s-1)(s-x)]^{1/2}} \frac{1}{((\frac{x-y}{[(s-1)(s-x)]^{1/2}})^2+1)^{3/2}} f(y) (1-y)^\alpha dy \, ds\\
&=&  (1-r) \int_k^{2} \frac{1}{(s-k)^{1/2}(s-1)(s-x)}  \int_0^x (s-y)^{-\alpha}(s-y)\\
&&\quad \quad \quad \quad  \quad \quad  \quad  \times  \frac{1}{[(s-1)(s-x)]^{1/2}} \frac{1}{((\frac{x-y}{[(s-1)(s-x)]^{1/2}})^2+1)^{3/2}} f(y) (1-y)^\alpha dy \, ds\\
\end{eqnarray*}
Now writing 
$$ s-y = (s-x) + (x-y),$$
we get two terms. The first term, since $(s-y)^{-\alpha} \leq (s-x)^{-\alpha}$ is then the same as in i-1) i.e. we get the right bound in that case. For the second term, we get that is bounded by
\begin{eqnarray*}
&&(1-r) \int_k^{2} \frac{1}{(s-k)^{1/2}(s-1)^{1/2}(s-x)^{1/2}}  \int_0^x (s-y)^{-\alpha}\\
&&\quad \quad \quad \quad  \quad \quad  \quad  \times  \frac{1}{[(s-1)(s-x)]^{1/2}} \frac{1}{(\frac{x-y}{[(s-1)(s-x)]^{1/2}})^2+1} f(y) (1-y)^\alpha dy \, ds\\
&=&(1-r) \int_k^{2} \frac{1}{(s-k)^{1/2}(s-1)^{1/2}(s-x)^{1/2}}  \int_0^x (s-y)^{-\alpha}\\
&&\quad \quad \quad \quad  \quad \quad  \quad  \times  \frac{1}{[(s-1)(s-x)]^{1/2}} \frac{1}{(\frac{x-y}{[(s-1)(s-x)]^{1/2}})^2+1} f(y) (1-y)^\alpha dy \, ds\\
\end{eqnarray*}
Then, 
$$  \frac{1}{[(s-1)(s-x)]^{1/2}} \frac{1}{(\frac{x-y}{[(s-1)(s-x)]^{1/2}})^2+1}  =  \frac{1}{\lambda} k_4(\frac{x-y}{\lambda}) $$
been a Poisson type kernel, the expression in the inner integral satisfies a (unilateral) condition of  Zygmund's lemma with respect to the measure $\mu(dy) = (1-y)^\alpha dy$ and therefore
\begin{eqnarray*}
&&J_\alpha f(x)  \leq C f^*_{J^{\alpha,\beta}}(x) (1-r) \int_k^{2} \frac{1}{(s-k)^{1/2}(s-1)}  \int_0^x  \frac{1}{(s-x)^\alpha} \frac{1}{\lambda} k_3(\frac{x-y}{\lambda}) (1-y)^\alpha dy \, ds,
\end{eqnarray*}
where $\lambda= [(s-1)(s-x)]^{1/2}$. Then by i-2) of the proof of  Theorem \ref{mainest}   gives us that the last term bounded, i. e.
$$ J_\alpha f(x)  \leq C f^*_{J^{\alpha,\beta}}(x).$$ 
\end{enumerate}
\item[ii)]  Case $-1< \alpha< 0.$
\begin{enumerate}
\item[ ii-1)] If $y\geq x$: In this range we have,
\begin{eqnarray*}
&&J_\alpha f(x)  \leq (1-r) \int_k^{2} \frac{1}{(s-k)^{1/2}(s-1)}  \int_x^1  \frac{1}{(s-x)^\alpha}\\
&&\quad \quad \quad \quad  \quad \quad  \quad  \times  \frac{1}{[(s-1)(s-x)]^{1/2}} \frac{1}{((\frac{x-y}{[(s-1)(s-x)]^{1/2}})^2+1)^{3/2}} f(y) (1-y)^\alpha dy \, ds
\end{eqnarray*}
Then, as in the case i-1), using the kernel $k_2$ by Zygmund's lemma with respect to the measure $\mu(dy) = (1-y)^\alpha dy$ \begin{eqnarray*}
&&J_\alpha f(x)  \leq C f^*_{J^{\alpha,\beta}}(x) (1-r) \int_k^{2} \frac{1}{(s-k)^{1/2}(s-1)}  \int_x^1  \frac{1}{(s-x)^\alpha} \frac{1}{\lambda} k_2(\frac{x-y}{\lambda}) (1-y)^\alpha dy \, ds,
\end{eqnarray*}
where $\lambda= [(s-1)(s-x)]^{1/2}$. Then by   ii-1)  of the proof of  Theorem \ref{mainest}, we get the last term is bounded i.e.
$$ J_\alpha f(x)  \leq C f^*_{J^{\alpha,\beta}}(x).$$
\item[ ii-2)] If $0< y< x$: In this range we have,
\begin{eqnarray*}
&&J_\alpha f(x)  \leq (1-r) \int_k^{2} \frac{1}{(s-k)^{1/2}(s-1)(s-x)}  \int_0^x  (s-y)^{1-\alpha}\\
&&\quad \quad \quad \quad  \quad \quad  \quad  \times  \frac{1}{[(s-1)(s-x)]^{1/2}} \frac{1}{((\frac{x-y}{[(s-1)(s-x)]^{1/2}})^2+1)^{3/2}} f(y) (1-y)^\alpha dy \, ds.
\end{eqnarray*}
Now since $\alpha <0$
$$(s-y)^{1-\alpha} \leq C_\alpha [(s-x)^{1-\alpha}+(x-y)^{(1-\alpha)}],$$
 we get two terms
\begin{eqnarray*}
&&C_\alpha  (1-r) \int_k^{2} \frac{(s-x)^{-\alpha}}{(s-k)^{1/2}(s-1)}\\
&&\quad \quad \quad \quad \times \int_0^x  \frac{1}{[(s-1)(s-x)]^{1/2}}\frac{1}{((\frac{x-y}{[(s-1)(s-x)]^{1/2}})^2+1)^{3/2}} f(y) (1-y)^{\alpha} dy \, ds \\
& &\quad \quad \quad + C_\alpha (1-r) \int_k^{2} \frac{1}{(s-k)^{1/2}(s-1)(s-x)}\\
&&\quad \quad \quad \quad \times \int_0^x  \frac{(x-y)^{1-\alpha}}{(s-1)^{1/2}(s-x)^{1/2}((\frac{x-y}{[(s-1)(s-x)]^{1/2}})^2+1)^{3/2}}  f(y) (1-y)^{\alpha}dy \, ds
\end{eqnarray*}
The first integral can be handle  in a similar way as in ii-1) using the kernel $k_2$ and ii-2) of the proof of  Theorem \ref{mainest}. For the second integral, by the similar argument as in ii-2) of the proof of Theorem \ref{mainest}, we have the bound,
\begin{eqnarray*}
&& (1-r) \int_k^{2} \frac{(s-1)^{-\alpha/2} (s-x)^{-\alpha/2}}{(s-k)^{1/2}}\\
&&\quad \quad \times \int_0^x  \frac{1}{(s-1)^{1/2}(s-x)^{1/2}}\frac{1}{[(\frac{x-y}{[(s-1)(s-x)]^{1/2}})^2+1]^{3/2-(1-\alpha)/2}}f(y) (1-y)^{\alpha}dy \, ds\\
\end{eqnarray*}
Now, considering the Poisson type kernel $k_4(x) = \frac{1}{(x^2+1)^{3/2-(1-\alpha)/2}},$ the expression in the inner integral satisfies a (unilateral) condition of  Zygmund's lemma with respect to the measure $\mu(dy) = (1-y)^\alpha dy$ and therefore
\begin{eqnarray*}
&&J_\alpha f(x)  \leq C f^*_{J^{\alpha,\beta}}(x) (1-r) \int_k^{2} \frac{1}{(s-k)^{1/2}(s-1)}  \int_0^x  \frac{1}{(s-x)^\alpha} \frac{1}{\lambda} k_2(\frac{x-y}{\lambda}) (1-y)^\alpha dy \, ds,
\end{eqnarray*}
where $\lambda= [(s-1)(s-x)]^{1/2}$. Then by   ii-2)  of the proof of  Theorem \ref{mainest} gives us that the last term bounded, i. e.
$$ J_\alpha f(x)  \leq C f^*_{J^{\alpha,\beta}}(x).$$ 
 \end{enumerate}
\end{enumerate}
\ep

{\bf Observation} Observe that there is an analogous operator 
 \begin{eqnarray}
\nonumber J_\beta f(x) &=&\int_{-1}^0 (1-r) \int_k^{2} \frac{(s-\min(x,y))^{1-\alpha}}{((x-y)^2+(s-1)(s- \min(x,y)))^{3/2}} \frac{ds}{(s-k)^{1/2}} \\
&& \quad \quad \quad  \quad \quad \quad \quad \quad  \quad \quad \quad  \quad \quad \quad \quad  \quad \quad \quad  \times(1+y)^{\beta} f(y) \, dy,
\end{eqnarray}

With analogous arguments as  in the previous result we have immediately 
 \begin{equation}
J_\beta f(x) \leq  C f^*_{J^{\alpha,\beta}}(x).
\end{equation}
Therefore, by the continuity properties of $f^*_{J^{\alpha,\beta}}$,  we have
\begin{corollary} The operators
$J_{\alpha}$ and $ J_\beta$ are  weak-$(1,1)$ continuous with respect to the Jacobi measure $J^{\alpha,\beta}$
\end{corollary}
Then, using the inequality (\ref{basicineq}) and the two previous results, we get that Jacobi maximal function $f^*_{\alpha, \beta}$ (see \ref{JacobMaxfunt}) is weak $(1,1)$ with respect to the Jacobi measure.\\

Now, let us consider  a Calder\'on-Zygmund's decomposition for a non-atomic Borel measure $\mu$ on $\mathbb{R}$ 
\begin{theorem} (Calder\'on-Zygmund) Given   $-\infty \leq a < b \leq \infty$, a non-atomic Borel measure $\mu$ with support on $(a,b)$, $\lambda >0$ and $f \in L^1(\mu)$, $f \geq0$, then   there exists a family of non-overlapping intervals $\{I_k\}$
\begin{enumerate}
\item [i)] $ \lambda < \frac{1}{\mu(I_k)} \int_{I_k} f(y)\mu(dy) \leq 2 \lambda,$
\item [ii)] $|f(x)| \leq \lambda, \, a.e. \mu,$ for $ x \notin \cup_k I_k.$
\end{enumerate}
\end{theorem}

\dem
\begin{itemize}

\item If $\frac{1}{\mu(a,b)} \int_a^b f(y) \mu(dy) >\lambda$ then $$\mu(a,b)<\frac{1}{\lambda} \int_a^b f(y) \mu(dy) = \frac{1}{\lambda} \| f\|_1,$$
and then there is nothing to prove.\\

\item  If $\frac{1}{\mu(a,b)} \int_a^b f(y) \mu(dy) \leq \lambda$ then consider two  intervals, $I_{0,1}, \, I_{0,2}$ with disjoint interiors such that $ (a,b) = I_{0,1} \cup I_{0,2}$ and $\mu(I_{0,1})=\mu(I_{0,2})=\frac{1}{2} \mu(a,b)$. Let us observe that we can not have that the inequality
$$ \frac{1}{\mu(I_{0,i})} \int_{ I_{0,i}} f(y) \mu(dy)> \lambda,$$
hold for both $i=1$ and $i=2$ since otherwise, 
\begin{eqnarray*}
 \frac{1}{\mu(a,b)} \int_{(a,b) } f(y) \mu(dy)=  \frac{2}{\mu(I_{0,1})} \int_{ I_{0,1}} f(y) \mu(dy)+  \frac{2}{\mu(I_{0,2})} \int_{ I_{0,2}} f(y) \mu(dy)> 4\lambda.
\end{eqnarray*}
which is a contradiction, then we have that at least one of then  (or even both)  satisfy
$$ \frac{1}{\mu(I_{0,i})} \int_{ I_{0,i}} f(y) \mu(dy)\leq \lambda.$$
In that case consider again two  intervals, $I_{i,1}, \, I_{i,2}$ with disjoint interiors such that $  I_{0,i} = I_{i,1} \cup I_{i,2}$ and $\mu(I_{i,1})=\mu(I_{i,2})=\frac{1}{2} \mu( I_{0,i}) = \frac{1}{4} \mu(a,b)$ and iterate the previous argument. If we have
$$ \frac{1}{\mu(I_{0,i})} \int_{ I_{0,i}} f(y) \mu(dy)> \lambda,$$
then
\begin{eqnarray*}
\frac{1}{\mu(I_{0,i})} \int_{ I_{0,i}} f(y) \mu(dy) &\leq& \frac{1}{\mu(I_{0,i})} \int_{(a,b)} f(y) \mu(dy)\\
&=&   \frac{2}{\mu(a,b)} \int_{ (a,b)} f(y) \mu(dy)\leq 2\lambda.
\end{eqnarray*}
Set $I_{0,i} $ aside, it will be one  of our chosen interval $I_k$. \\
This infinite recursion will give us a family $\{I_k\}$ such that,
$$ \lambda < \frac{1}{\mu(I_k)} \int_{I_k} f(y)\mu(dy) \leq 2 \lambda.$$

Set $G_\lambda = \cup_{k=1}^\infty I_k$, then
\begin{eqnarray*}
\mu(G_\lambda) &=& \sum_{k=1}^\infty \mu(I_k) < \frac{1}{\lambda} \sum_{k=1}^\infty \int_{I_k} f(y)\mu(dy)\leq  \frac{1}{\lambda} \int_{G_\lambda} f(y)\mu(dy)\\
 &\leq& \frac{1}{\lambda} \int_{\mathbb{R}} f(y)\mu(dy) =  \frac{1}{\lambda} \|f\|_{1,\mu}.
\end{eqnarray*}
Let us observe that if  $x \notin  \cup_k I_k$ then there is an infinite family of intervals $I$ containing $x$ such that 
$$ \frac{1}{\mu(I)} \int_{ I} f(y) \mu(dy)\leq \lambda,$$
then by Lebesgue differentation theorem,  see Lemma 7 of \cite{apCal}, we get  $|f(x)| \leq \lambda$ a.e.$ \mu$, $x \notin  \cup_k I_k$. \\

Now, set  $\mu_k = \frac{1}{\mu(I_k)}\int_{I_k} f(y)\mu(dy)$ we can write  $f = g+ b$ where, 
$$g(x)=f \chi_{\mathbb{R} -G_\lambda} (x)+  \sum_k \mu_k \chi_{I_k}(x)$$
and 
$$ b (x) = f(x) - g(x) =  \sum_k  (f(x) - \mu_k) \chi_{I_k}(x).$$
 $g, b $ are called  that good and bad part of $f$ respectively. Observe that $ g \leq 2 \lambda$ in $G_\lambda$, the bad part is only non-zero in $G_\lambda$ and $\int_{I_k} b(y) \mu(dy) = 0$.\\

\end{itemize}
If  $G^*_\lambda = \cup_{k=1}^\infty I^*_k$ where $I^*_k = 3 I_k$ meaning that $I^*_k$ is the union of $I_k$ with two other intervals (one to the right and one to the left of it) with the same $\mu$ measure, i.e. $I^*_k = I'_k \cup I_k \cup I''_k$, with $\mu(I'_k)= \mu (I_k) = \mu (I''_k),$ then
$$ \mu(G^*_\lambda) = \sum_{k=1}^\infty \mu(I^*_k) =3 \sum_{k=1}^\infty \mu(I_k)\leq  \frac{3}{\lambda} \|f\|_{1,\mu}.$$ \ep

We can use Calder\'on-Zygmund decomposition  for a kernel $K(r, x,y)$  that satisfies the conditions of Zygmund's lemma 

\begin{proposition} Given  a non-atomic Borel measure $\mu$, with support in $(a,b)$, and a kernel $K(r, x,y)$  that satisfies  Zygmund's lemma conditions  (\ref{Zygbound1}) and (\ref{Zygbound2}) with respect to $\mu$, i. e. 
$$
 \int_a^b |K(r,x,y)| \mu(dy) \leq M_1
$$
and 
$$
 \int_x^b \mu(x,y) V_2(K(r, x,dy)) \leq M_2, \quad \int_a^x  \mu(y,x) V_2(K(r, x,dy)) \leq M_2.
$$
Then  for $f \in L^1(\mu)$ and $ x \notin G^*_\lambda$,
\begin{equation}
 \sup_r |  \int_a^b K(r, x,y) f(y) \mu(dy)| \leq C \lambda.
\end{equation}
\end{proposition} 

\dem

We know by Zygmund's lemma that,
$$|\int_a^b K(r, x,y) f(y) \mu(dy)| \leq M f_\mu^*(x).$$
Now, using Calder\'on-Zygmund decomposition for $f = g+b$, we get
$$ \int_a^b K(r, x,y) f(y) \mu(dy) = \int_a^b K(r, x,y) g(y) \mu(dy) + \int_a^b K(r, x,y) b(y) \mu(dy)$$
and as $|g| < 2 \lambda$, a.e.$ \mu$
$$ |\int_a^b K(r, x,y) g(y) \mu(dy)| < 2 M_1 \lambda.$$
If $ x \notin G^*_\lambda$ using integration by parts, where $I_k =(a_k, b_k)$ 
\begin{eqnarray*}
| \int_a^b K(r, x,y) b(y) \mu(dy) | &=& |\sum_k  \int_{I_k} (f(y) - \mu_k) K(r, x,y)  \mu(dy) | \\
&=& |\sum_k  \int_{a_k}^{b_k} (f(y) - \mu_k) K(r, x,y)  \mu(dy) | \\
&=& |\sum_k   \int_{a_k}^{b_k}( \int_{a_k}^{y} (f(u) - \mu_k) \mu(du) ) K(r, x,dy)  | \\
\end{eqnarray*}
as $ \int_{a_k}^{b_k} (f(u) - \mu_k) \mu(du) = 0,$  and using that $ x \notin G^*_\lambda$,
\begin{eqnarray*}
| \int_a^b K(r, x,y) b(y) \mu(dy) |& \leq & C \lambda \sum_k \int_{I_k}  \mu(I_k) V_2(K(r, x,dy)) \\
&\leq & C \lambda \sum_k \int_{I_k}  \mu(x,y ) V_2(K(r, x,dy)) \\
&\leq & C \lambda \int_{G_\lambda}  \mu(x,y ) V_2(K(r, x,dy)) \\
&\leq & C \lambda \int_x^b  \mu(x,y ) V_2(K(r, x,dy)) \leq C \lambda M_2. \\
\end{eqnarray*}
Thus, for $ x \notin G^*_\lambda$
$$ \sup_r |  \int_a^b K(r, x,y) f(y) \mu(dy)| \leq C \lambda.$$
\ep \\
This result could be extended to the case of measures that do have atoms.\\

The following result was proved implicitly by L. Cafarelli in \cite{cafacal}, 

\begin{theorem}
The Jacobi measure $J^{\alpha,\beta}$ is a doubling measure.
\end{theorem}

\dem Let us consider first the measure $\mu(dy) = y^a,$  in $[0,1], \,  a >-1$. Then we will see that $\mu$ is a doubling measure on $[0,1]$. 

Let $k \geq 2$ and $I_{k,j} = [ k 2^{-j}, (k+1) 2^{-j}]$ a dyadic interval.
Observe that
$$\mu(I_{k,j})= \int_{k 2^{-j}}^{(k+1) 2^{-j}} y^a dy = \frac{2^{-j(a+1)}}{a+1}[(k+1)^{a+1}- k^{a+1}]$$
Now let us consider 3$I_{k,j}$ the interval with the same center $(k+1/2) 2^{-j}$ and 3 times the length 
of  $I_{k,j}$ i.e. $3 I_{k,j} = [ (k-1) 2^{-j}, (k+2) 2^{-j}]$, then
$$\mu(3I_{k,j})= \int_{(k-1) 2^{-j}}^{(k+2) 2^{-j}} y^a dy = \frac{2^{-j(a+1)}}{a+1}[(k+2)^{a+1}- (k-1)^{a+1}].$$
Thus,
\begin{eqnarray*}
\frac{\mu(3I_{k,j})}{\mu(I_{k,j})} &=& \frac{(k+2)^{a+1}- (k-1)^{a+1}}{(k+1)^{a+1}- k^{a+1}} \\
&=& \frac{(1+\frac{2}{k})^{a+1}- (1-\frac{1}{k})^{a+1}}{(1+\frac{1}{k})^{a+1}-1}. \\
\end{eqnarray*}
It can be proved that the quotient $\frac{\mu(3I_{k,j})}{\mu(I_{k,j})} $ is increasing in $k$ for $a \in (0,1)$ and decreasing for $a \in (-1,0) \cup (1,\infty)$. By L'Hopital rule,
\begin{eqnarray*}
\lim_{k \rightarrow \infty} \frac{\mu(3I_{k,j})}{\mu(I_{k,j})} &=&\lim_{k \rightarrow \infty}  \frac{2(1+\frac{2}{k})^{a}+ (1-\frac{1}{k})^{a}}{(1+\frac{1}{k})^{a}} = 3.
\end{eqnarray*}
Therefore if $ a \in (0,1)$, 
$$C_a=\frac{3^{a+1}}{2^{a+1} -1} \leq  \frac{\mu(3I_{k,j})}{\mu(I_{k,j})}  \leq 3.$$
and elsewhere
$$3\leq \frac{\mu(3I_{k,j})}{\mu(I_{k,j})}  \leq \frac{3^{a+1}}{2^{a+1} -1}= C_a.$$

Similarly, using the same arguments we can prove that  $\mu$ is also a doubling measure on $[-1,0]$. 

Now observe that, by a change of variable,  on $[0,1]$ the measure $y^a dy$ is equivalent to $(1-y)^a dy$, in the following sense
$$ \int_0^1 f(y) (1-y)^\alpha dy = \int_0^1 f(1-u) u^\alpha du = \int_0^1\overline{f(u)} u^\alpha du,$$
and clearly there is a one-to-one correspondence between $f$ and $\overline{f}$.\\
 Similarly, on $[-1,0]$ the measure $y^a dy$ is equivalent to $(1+y)^a dy,$

Finally, as a consequence of the previous results  we have that  the Jacobi measure $J^{\alpha, \beta}(dy) = (1-y)^{\alpha}(1+y)^{\beta} dy$ in $(0,1)$  is equivalent to $ y^\alpha dy$ and  is equivalent to $ y^\beta dy$ in $(-1,0)$. Therefore $J^{\alpha, \beta}$ is then a doubling measure on $[-1,1]$. \ep \\

Now that we know that the Jacobi measure $J^{\alpha,\beta}$ is a doubling measure we can use the result of A. P. Calder\'on \cite{apCal},  in order to get the $A_p$ weight theory for $J^{\alpha,\beta}$. Remember a  function $\omega >0$, is an $A_p$ weight, $\omega \in A_p$,  if

\begin{equation}
[\frac{1}{J^{\alpha,\beta}(B)}\int_B \omega(y) J^{\alpha,\beta}(dy)][ \frac{1}{J^{\alpha,\beta}(B)}\int_B \omega(y)^{-1/(p-1)} J^{\alpha,\beta}(dy)]^{p-1} \leq C_p, 
\end{equation}
for $1< p < \infty$ and 

\begin{equation}
\omega^*_{J^{\alpha,\beta}} (x) \leq C_1 \omega(x),
\end{equation}
for $p=1.$

For a complete exposition of the $A_p$ weight theory see for instance  Duoandikoetxea \cite{duo}.\\
In what follows we will use the following notation for a measure $\mu(dx) = g(x) dx$,
$$ \int_a^b \mu(dy) = \int_a^b g(y) dy = G(b) - G(a).$$

We want to consider some interesting $A_1$ weights for the Jacobi measure. Observe that by the factorization result (see  Duoandikoetxea \cite{duo}, Proposition 7.2, page 136) they are like building blocks for $A_p$ weights for $ p >1$.  First of all, we need the following technical result. 
\begin{lemma}\label{maxfunEval}
Let  $\mu$ be a non-negative Borel measure on $[0,1)$ and absolutely  continuos i.e. $\mu(dx) = g(x) dx$ where $g$ is non-negative and continuous. Then if $f$ is a non increasing non negative function then 
$$ \frac{1}{G(x)- G(a)} \int_x^a f(y) g(y) dy$$
 is also non-increasing function. \\
The same result is true for a  non-negative Borel  measure $\mu$ on $(-1,0]$.
 \end{lemma}
 
 \dem
 Since 
 \begin{eqnarray*}
\frac{d}{dx}\big( \frac{1}{G(a)- G(x)} \int_x^a f(y) g(y) dy \big) &=& \frac{-f(x)g(x) (G(a) - G(x))+ (\int_x^a f(y) g(y) dy) g(x)}{(G(a) -G(x))^2}\\
&=& \frac{g(x)(-f(x) (G(a) - G(x))+ \int_x^a f(y) g(y) dy)}{(G(a) -G(x))^2} \leq 0,\\
\end{eqnarray*}
as $ g \geq 0$ and $ f(x)  \int_x^a g(y) dy \geq  \int_x^a f(y) g(y) dy.$

 Therefore the quotient is non-increasing
as claimed. \ep \\

We will use the previous result to consider lateral maximal functions. If we consider the left lateral maximal function of  non increasing non negative function $f$,
$$f^*_{-} (a) = \sup_{0 \leq x\leq a} \frac{1}{G(a)- G(x)} \int_x^a f(y) g(y) dy,$$
we have, by Lemma \ref{maxfunEval},
$$ f^*_{-} (a)  =  \frac{1}{G(a)- G(0)} \int_0^a f(y) g(y) dy =\frac{1}{G(a)} \int_0^a f(y) g(y) dy,$$
as $G(0) =0$.

By analogous argument we have that for a non increasing non negative function $f$, its the right lateral maximal function equals,
$$ f^*_{+} (a)  =  \sup_{0 \leq a\leq x} \frac{1}{G(x)- G(a)} \int_a^x f(y) g(y) dy = f(a_+).$$

The case of a general non-negative Borel  measure $\mu$ can be obtained using Helly's selection principle.\\

Let us finally consider the $A_1$  weights for the Jacobi measure,

\begin{lemma}\label{jacobiweight}

i) For $1 < \alpha < \infty$, let us consider the power measure $\mu_\alpha(dx)=x^\alpha dx$ on $[0,1)$, then  the measure $\mu_{\overline{\alpha}}(dx) = x^{\overline{\alpha}} dx, \, -1 < \overline{\alpha} < 0, \, \alpha + \overline{\alpha} >-1$ is a $A_1$ weight with respect to $\mu_\alpha$. \\
ii) Similarly, considering the power measure $\mu_\beta(dx)=x^\beta$ on $[-1,0)$, then  the measure $\mu_{\overline{\beta}}(dx) =x^{\overline{\beta}} dx, \, -1 < \overline{\beta} < 0, \, \beta + \overline{\beta} >-1$  is an $A_1$ weight with respect to the $\mu_\beta$.
\end{lemma}

\dem
By previous considerations, the left  maximal function with respect to $\mu_\alpha$ is equal to,
$$ \frac{C}{x^{\alpha +1}} \int_0^x t^\alpha t^{\overline{\alpha}} dt =  \frac{C}{x^{\alpha +1}} \int_0^x t^{\alpha+\overline{\alpha}} dt = \frac{C}{x^{\alpha +1}} x^{\alpha+\overline{\alpha}+1} = C x^{\overline{\alpha}},$$
and from the right is simply $x^{\overline{\alpha}}$,
i.e.  the measure $\mu_{\overline{\alpha}}(dx)= x^{\overline{\alpha}} dx, \, -1 < \overline{\alpha} < 0, \, \alpha + \overline{\alpha} >-1$  is an $A_1$ weight with respect to the measure $\mu_\alpha$.

Similarly, on $[-1,0)$ $\mu_{\overline{\beta}}(dx) = x^{\overline{\beta}} dx, \, -1 < \overline{\beta} < 0, \, \beta + \overline{\beta} >-1$  is an $A_1$ weight with respect to the measure $\mu_\beta(dx) =x^\beta dx,\, 1 < \beta < \infty $. \ep\\
Now we have the following result for the Jacobi measure.  This result extends the set of weights that were considered in \cite{calv}, where only positive power were considered.
\begin{theorem}\label{Jacobiweight} The measures 
\begin{equation}
\mu_{\overline{\alpha},\overline{\beta} }(dx)= (1-x)^{\overline{\alpha}}(1+x)^{\overline{\beta}} dx, \; \overline{\alpha} +\overline{\beta} >-1,
\end{equation}
are $A_1$ weights with respect to the Jacobi measure  $J^{\alpha,\beta}(dx) = (1-x)^\alpha(1+x)^\beta dx, \, \alpha +\beta >-1$
on $[-1,1]$ 
\end{theorem}
\dem
By Lemma \ref{jacobiweight} and similar arguments as above, the measure $\nu_{ \overline{\alpha}}(dx) =(1-x)^{ \overline{\alpha}}dx, \, -1 < \overline{\alpha} < 0, \;  \alpha+ \overline{\alpha} >-1$  is an $A_1$ weight with respect to the measure $\nu_\alpha(dx) = (1-x)^\alpha dx$ and similarly, the measure $\nu_{ \overline{\beta}}(dx) =(1+x)^{ \overline{\beta}}dx, \; -1 < \overline{\beta} < 0, \,  \beta + \overline{\beta} >-1 >-1$  is an $A_1$ weight with respect to the measure $\nu_\beta(dx) = (1+x)^\beta dx$ on $[-1,0)$ and from there we get our result inmediately. \ep \\

Finally, as a corollary of Theorem \ref{Jacobiweight} we have the following result for Abel summability of Jacobi function expansions. 
\begin{corollary} For the Abel summability of Jacobi function expansions we have for $1 < p < \infty$,
$$||  \tilde{f}^{\alpha, \beta}(r, \cdot) ||_p \leq C || f ||_p$$
\end{corollary}

\dem

Let us consider only the case of the interval $[0,1]$ with $\alpha<0$, the case $[-1,0]$ is totally analogous. From (\ref{JabFuncIntRep2}) we have, by the maximal inequality of the Hardy-Littlewood function $M_{\nu_\alpha} $ with respect to the measure $\nu_\alpha(dx) = (1-x)^\alpha dx $,
\begin{eqnarray*}
 \int_0^1\tilde{K}^{\alpha,\beta}(r,x,y) f(y) dy &\leq & C_\beta (1-x)^{\alpha/2} \int_0^1 K^{\alpha,\beta}(r,x,y)(1-y)^{\alpha/2}  f(y) dy \\
 &=& C_\beta(1-x)^{\alpha/2} \int_0^1 K^{\alpha,\beta}(r,x,y) [(1-y)^{-\alpha/2} f(y)] (1-y)^{\alpha}dy\\
 &\leq& C_\beta M_{\nu_\alpha}  ((1-\cdot)^{-\alpha/2} f) (x) (1-x)^{\alpha/2}.
\end{eqnarray*}
Therefore, by the $L^2$ continuity of  $M_{\nu_\alpha} $  with respect to the measure $\nu_\alpha$,
\begin{eqnarray*}
 \int_0^1[ \int_0^1\tilde{K}^{\alpha,\beta}(r,x,y) f(y) dy ]^2 dx
 &\leq& C_\beta \int_0^1 [M_{\nu_\alpha}  ((1-\cdot)^{-\alpha/2} f)]^2 (x) (1-x)^{\alpha} dx\\
&\leq & C_\beta \int_0^1  [(1-x)^{-\alpha/2} f(x) ]^2 (1-x)^\alpha dx\\ 
&\leq& C_\beta  \int_0^1  [f(y)]^2  dy= C \| f\|_2^2.
\end{eqnarray*}
Thus
$$ ||  \tilde{f}^{\alpha, \beta}(r, \cdot) ||_2 \leq C || f ||_2.$$

Analogously, for the $L^p$ inequality. If $p>2$,  
\begin{eqnarray*}
\int_0^1[ \int_0^1\tilde{K}^{\alpha,\beta}(r,x,y) f(y) dy ]^p dx
 &\leq& C_\beta \int_0^1[M_{\nu_\alpha} ( (1-\cdot)^{-\alpha/2} f)]^p(x)  (1-x)^{p\alpha/2}dx\\
 \end{eqnarray*}
 and observe that
 $$ (1-x)^{p\alpha/2} = (1-x)^{p\alpha/2 -\alpha+\alpha} = (1-x)^{\delta + \alpha}, \quad \delta = p\alpha/2 -\alpha= \delta/2 (p-2);$$
 $(1-x)^{\delta + \alpha}$ is a $A_p (\nu_\alpha )$-weight  if $\delta + \alpha = \alpha p/2 >-1,$ i.e. $p < 2/|\alpha|$, and therefore, by the $L^p$ continuity of  $M_{\nu_\alpha}  $ with respect to the measure $\nu_\alpha$,
 \begin{eqnarray*}
 \int_0^1[ \int_0^1\tilde{K}^{\alpha,\beta}(r,x,y) f(y) dy ]^p dx
&\leq& C_\beta  \int_0^1 (1-x)^{-\alpha p/2} [f(x)]^p (1-x)^{\alpha p/2}  dx\\
&=& C_\beta  \int_0^1 [ f(x)]^p  dx= C_\beta  \| f\|_p^p.
\end{eqnarray*}

If $1 < p <2$, $ (1-x)^{p\alpha/2 -\alpha}$ is a $A_p (\nu_\alpha )$-weight if and only if 
$$ (1-x) ^{(p\alpha/2 -\alpha)(-1/(p-1))} = (1-x) ^{(-\alpha/2 p/(p-1)+\alpha/(p-1)}  =(1-x) ^{(-\alpha/2 q+\alpha(q-1))} ,$$
is a $A_q (\nu_\alpha)$-weight, $\frac{1}{p}+\frac{1}{q}=1$, see \cite{duo}. But
$$ -\alpha/2 q+\alpha(q-1)= -\alpha/2 q+\alpha q- \alpha =  \alpha/2 q- \alpha=\gamma,$$
and therefore $ (1-x) ^{(p\alpha/2 -\alpha)(-1/(p-1))} = (1-x) ^\gamma$ is a $A_q (\nu_\alpha dy)$-weight for $q >2$. Then
 \begin{eqnarray*}
 \int_0^1[ \int_0^1\tilde{K}^{\alpha,\beta}(r,x,y) f(y) dy ]^q dx
&\leq& C_\beta  \int_0^1 (1-x)^{-\alpha q/2} [f(x)]^q (1-x)^{\alpha q/2}  dx\\
&=& C_\beta  \int_0^1 [ f(x)]^q  dx= C_\beta  \| f\|_q^q.
\end{eqnarray*}
From the previous case the condition $p < 2/|\alpha|$ holds if and only if $ q > \frac{2}{2-|\alpha|},$
so the general condition for $p$ is 
$$ \frac{2}{2-|\alpha|} < p < \frac{2}{|\alpha|}.$$

Therefore, the bilateral condition in $[-1,1]$ is  
$$\max[ \frac{2}{2-|\alpha|},\frac{2}{2-|\beta|}] < p < \min[ \frac{2}{|\alpha|}, \frac{2}{|\beta|}].$$\ep

\end{document}